\documentclass[11pt, final]{amsart}
\usepackage{amssymb}
\usepackage{showkeys}
\usepackage[mathscr]{eucal}

\theoremstyle{plain}
\newtheorem{theorem}{Theorem}
\newtheorem{lemma}{Lemma}

\theoremstyle{remark}
\newtheorem{remark}{Remark}

\def\cA{\mathscr{A}}
\def\cB{\mathscr{B}}
\def\C{\mathbb C}
\def\l{\lambda}
\def\la{\langle}
\def\N{\mathbb N}
\def\ot{\otimes}
\def\R{\mathbb R}
\def\ra{\rangle}
\def\rank{\operatorname{rank}}
\def\rng{\operatorname{rng}}
\def\tr{\operatorname{tr}}
\def\wlim{\operatorname{w-\lim}}

\begin{document}
\title[]
{Some Multiplicative Preservers on $B(H)$}
\author{LAJOS MOLN\'AR}
\address{Institute of Mathematics\\
         Lajos Kossuth University\\
         4010 Debrecen, P.O.Box 12, Hungary}
\email{molnarl@math.klte.hu}
\thanks{  This research was supported from the following sources:\\
          1) Joint Hungarian-Slovene research project supported
          by OMFB in Hungary and the Ministry of Science and
          Technology in Slovenia, Reg. No. SLO-2/96,\\
          2) Hungarian National Foundation for Scientific Research
          (OTKA), Grant No. F--019322,\\
          3) A grant from the Ministry of Education, Hungary, Reg.
          No. FKFP 0304/1997}
\date{\today}
\begin{abstract}
In this paper we describe the form of those continuous multiplicative
maps on $B(H)$
($H$ being a separable complex Hilbert space of dimension not less than
3) which preserve the rank, or the corank. Furthermore, we characterize
those continuous *-semigroup endomorphisms of $B(H)$ which are spectrum
non-increasing.
\end{abstract}
\maketitle

\section{Introduction}

The study of linear preserver problems has a long history. In fact, it
is one of the most active research areas in matrix theory
\cite{LiTsing} (also see \cite{BrSe} for a survey on linear preservers
on operator algebras). In a recent paper \cite{Hoch} Hochwald started to
investigate multiplicative preserver problems. In his paper he described
the form of those multiplicative selfmaps of a matrix algebra which
preserve
the spectrum (also see \cite{Aup} for a result concerning Banach
algebras). As a natural generalization,
he also raised the question of spectrum-preserving multiplicative maps
on operator algebras even under the possible additional condition of
surjectivity. However, taking into account
Martindale's purely algebraic result \cite[First Corollary]{Mar}, it
follows that
in the case of many operator algebras $\cA$ (for example, if $\cA$ is a
standard operator algebra on a Banach space of dimension at least 2,
i.e. a subalgebra of the whole operator algebra which
contains the ideal of all finite rank operators), every
multiplicative transformation on $\cA$ which maps onto an arbitrary
algebra is automatically additive. Since additivity is not so far
from linearity, it seems a much more exciting
problem to try to attact the problem if surjectivity is not assumed.

In the present paper we consider those multiplicative preservers on
the operator algebra $B(H)$ which are the natural analogues of the most
'popular' linear preservers, that is, we try to describe the form
of those
multiplicative maps which preserve the rank, or the spectrum. Our main
tool on the way to obtain our results is the extensive theory of
measures on lattices of projections of operator algebras.

Let us fix the notation and the concepts that we shall use throughout.
Let $H$ be a Hilbert space. Denote by $B(H)$ the algebra of all bounded
linear operators on $H$. An operator $P\in B(H)$ is called an idempotent
if $P^2=P$. Two idempotents $P,Q\in B(H)$ are said to be orthogonal if
$PQ=QP=0$. We denote $P\leq Q$ if $PQ=QP=P$. Any self-adjoint idempotent
in $B(H)$ is called a projection.
The set of all projections in $B(H)$ is denoted by $P(H)$. The notation
$P_1(H)$ stands for the set of all rank-one projections on $H$.
The ideal of all finite-rank elements in $B(H)$ is denoted by $F(H)$.
If $x, y\in H$, then the operator $x \ot y$ is defined by
\[
(x\ot y)(z)=\la z,y\ra x \qquad (z \in H).
\]
Clearly, every rank-one operator $A$ is of the form $A=x \ot y$.
Moreover, the rank-one projections are exactly the operators of the form
$x\ot x$ where $x$ is a unit vector.

A linear map $\phi: \cA \to \cB$ between the algebras $\cA $ and $\cB$
is called a Jordan homomorphism if $\phi(x)^2=\phi(x^2)$ holds for
every $x\in \cA$.
Obviously, every homomorphism is a Jordan homomorphism and this is the
case with every antihomomorphism as well, that is, with every linear map
$\psi: \cA \to \cB$ for which $\psi(xy)=\psi(y)\psi(x)$ $(x,y \in \cA)$.

\section{Statement of the results}

In what follows let $H$ be a separable complex Hilbert space of
dimension at least 3.

Our first result describes the form of the continuous multiplicative
rank preservers
on $P(H)$. We emphasize that here we assign rank only to the
elements of $F(H)$.

\begin{theorem}\label{T:multrankPH}
Let $\phi: P(H) \to B(H)$ be a continuous multiplicative map which
preserves the rank. Then $\phi$ is of the form
\begin{equation}\label{E:formegy}
\phi(P)=TPS + \phi'(P) \qquad (P\in P(H))
\end{equation}
where $T,S:H \to H$ are either both bounded linear operators or both
bounded
conjugate-linear operators such that $ST=I$, $\phi' :P(H) \to B(H)$
is a continuous
multiplicative map which vanishes on the set of all finite rank
projections and $TPS\phi'(Q)=\phi'(Q)TPS=0$ $(P,Q\in P(H))$.
\end{theorem}

\begin{remark}
To see that the 'singular' part $\phi'$ can really appear in
\eqref{E:formegy} suppose that $H \cong H\oplus H$ (i.e., $H$ is
infinite dimensional) and define a continuous multiplicative map $\phi
:P(H) \to B(H)$ by
\[
\phi(P)=
\begin{cases}
\left[
\begin{matrix}
P & 0\\
0 & 0
\end{matrix}\right],
\quad \text{if the corank of $P$ is infinite}\\ \\
\left[
\begin{matrix}
P & 0\\
0 & I
\end{matrix}\right],
\quad \text{if the corank of $P$ is finite}.
\end{cases}
\]
\end{remark}

The description of continuous multiplicative rank preservers on $B(H)$
reads as follows.

\begin{theorem}\label{T:multrank}
Let $\phi: B(H) \to B(H)$ be a continuous multiplicative map which
preserves the rank. Then $\phi$ is of the form
\begin{equation}\label{E:formhar}
\phi(A)=TAS \qquad (A\in B(H))
\end{equation}
where $T,S: H\to H$ are either both bounded linear operators or both
bounded conjugate-linear operators and $ST=I$.
\end{theorem}

In the following result we describe the form of the continuous
multiplicative maps
on $B(H)$ which preserve the corank. There are (at least) two possible
definitions of the corank of an operator $A\in B(H)$. Let $n$ be a
nonnegative integer. The first possibility is as follows. We say that
the operator $A\in B(H)$ has corank $n$ if the algebraic dimension of
the quotient space $H/\rng A$ is $n$-dimensional. The second possibility
is when we say that $A$ has corank $n$ if the Hilbert dimension of
$\overline{\rng A}^\perp$ is $n$.
We shall see in the proof of our next result
that we have the same description in both cases.
In relation to the following theorem we also refer to \cite[Theorem
3]{GyMS} and \cite[Theorem 2]{Mol}.

\begin{theorem}\label{T:multcorank}
Let $\phi: B(H) \to B(H)$ be a continuous multiplicative map which
preserves the corank. Then $\phi$ is of the form
\begin{equation}\label{E:formnegy}
\phi(A)=TAT^{-1} \qquad (A\in B(H))
\end{equation}
where $T:H\to H$ is either a bounded linear operator or a bounded
conjugate-linear operator. In particular,
$\phi$ is either a linear or a conjugate-linear algebra automorphism
of $B(H)$.
\end{theorem}

Finally, we consider multiplicative maps $\phi$ on $B(H)$ that are
spectrum non-increasing which means that $\sigma (\phi(A)) \subset
\sigma (A)$ for every $A\in B(H)$.

\begin{theorem}\label{T:multspect}
Let $\phi: B(H) \to B(H)$ be a continuous *-semigroup homomorphism
(that is, a multiplicative map with the property that
$\phi(A)^*=\phi(A^*)$ $(A\in B(H))$). If $\phi$ is spectrum
non-increasing, then $\phi$
is a linear *-endomorphism of $B(H)$. More precisely, there are linear
isometries $U_n: H\to H$ $(n=1, \ldots)$ with pairwise orthogonal ranges
which generate $H$ such that $\phi$ is of the form
\begin{equation}\label{E:formot}
\phi(A)=\sum_n U_nAU_n^* \qquad (A\in B(H)).
\end{equation}
\end{theorem}

Notice that it is an obvious byproduct of the foregoing theorem that the
spectrum non-increasing maps appearing there are necessarily spectrum
preserving.

We remark that the form of our preservers in the low-dimensional cases
(i.e. when $\dim H\leq 2$) can be easily deduced from the result of \v
Semrl in \cite{Sem} where the general form of the multiplicative
selfmaps of a matrix algebra is given.

\section{Proofs}

\begin{proof}[Proof of Theorem~\ref{T:multrankPH}]
The idea of the proof is very simple. First extend $\phi$ to a linear
map on $F(H)$ (this will be denoted by $\tilde \psi$) which preserves
the rank and then apply a result on the form of linear rank preservers.
So, the idea is obvious but we have
to work hard to reach the desired conclusion.

Let $P_1, \ldots , P_n$ be pairwise orthogonal rank-one projections.
Let $P=P_1+ \cdots +P_n$. By the properties of $\phi$, $\phi(P_1),
\ldots, \phi(P_n)$
are pairwise orthogonal rank-one idempotents and $\phi(P)$ is a rank-$n$
idempotent. Since we have
\[
\begin{split}
\phi(P_1)+\cdots +\phi(P_n)= &
\phi(P)\phi(P_1)+\cdots +\phi(P)\phi(P_n)=\\
& \phi(P)(\phi(P_1)+\cdots +\phi(P_n))
\end{split}
\]
and
\[
\begin{split}
\phi(P_1)+\cdots +\phi(P_n)= &
\phi(P_1)\phi(P)+\cdots +\phi(P_n)\phi(P)= \\
& (\phi(P_1)+\cdots +\phi(P_n))\phi(P),
\end{split}
\]
it follows that $\phi(P_1)+\cdots +\phi(P_n)\leq \phi(P)$. But the
idempotents
on both sides of the latter inequality have rank $n$ which implies
that
\begin{equation}\label{E:hat}
\phi(P_1)+\cdots +\phi(P_n)= \phi(P).
\end{equation}
Let $H_d$ denote an arbitrary $d$-dimesional subspace of $H$.
Consider the natural embedding $B(H_d) \hookrightarrow B(H)$ and for any
$h\in H$ let $\phi_h$ be defined by $\phi_h(P)=\la \phi(P)h, h \ra$.
Taking \eqref{E:hat} into account we easily obtain
that $\phi_h$ is a measure on $P(H_d)$. We assert that $\phi$
is bounded on $P_1(H_d)$. Indeed, suppose on the contrary that there is
a sequence $(x_k)$ of unit vectors in $H_d$ such that $\| \phi(x_k \ot
x_k)\| \to \infty$. Since $H_d$ is finite
dimensional, $(x_k)$ has a convergent subsequence.
We can assume without any loss of generality that this subsequence is
the original sequence
$(x_k)$. Let $x=\lim_k x_k$. Then $x\in H_d$ is a unit vector and
by the continuity of $\phi$
we have $\| \phi(x \ot x)\| =\infty$ which is an obvious contradiction.
So, for any
$h\in H$, $\phi_h$ is a so-called $P_1$-bounded measure on $P(H_d)$. By
Gleason's theorem \cite[Theorem 3.2.16.]{Dvur} this implies that, in
case $d\geq 3$, there exists a linear operator $T_h$ on $H_d$ such that
\begin{equation}\label{E:gleas}
\phi_h (P)=\tr T_h P \qquad (P \in P(H_d)).
\end{equation}
Our aim now is to extend $\phi$ to a linear transformation of $F(H)$.
Let $x_1, \ldots , x_n\in H$ be unit vectors (the pairwise orthogonality
of the $x_i$'s is not assumed) and let $\l_1, \ldots , \l_n$ be real
numbers. Define \begin{equation}\label{E:het}
\psi(\sum_k \l_k x_k \ot x_k)=\sum_k \l_k \phi(x_k \ot x_k).
\end{equation}
We have to check that $\psi$ is well-defined. To see this, let
$y_1, \ldots, y_n\in H$ be unit vectors and
$\mu_1, \ldots , \mu_n \in \R$ such that
\[
\sum_k \l_k x_k \ot x_k=\sum_k \mu_k y_k \ot y_k.
\]
Let $H_d$ be a finite dimensional subspace of $H$ of dimension $d\geq
3$ which contains $x_1, \ldots , x_n , y_1, \ldots, y_n$. Let $h\in H$
be any vector. Let $T_h$ denote the linear operator on $H_d$
corresponding to $\phi_h$ (see \eqref{E:gleas}). We compute
\[
\la \sum_k \l_k \phi(x_k \ot x_k)h,h\ra =
\sum_k \l_k \phi_h(x_k \ot x_k) =
\tr T_h (\sum_k \l_k x_k \ot x_k)=
\]
\[
\tr T_h (\sum_k \mu_k y_k \ot y_k)=
\sum_k \mu_k \phi_h(y_k \ot y_k) =
\la \sum_k \mu_k \phi(y_k \ot y_k)h,h\ra .
\]
Since this holds true for every $h\in H$, we obtain that $\psi$ is
well-defined. The definition \eqref{E:het} now clearly implies that
$\psi$ is a real-linear operator on the set of all self-adjoint
finite-rank operators.
Clearly, $\psi$ sends projections to idempotents.
It is now a standard argument to verify that
the extension $\tilde \psi: F(H) \to
B(H)$ of $\psi$ defined by
\[
\tilde \psi(A+iB) =\psi(A) +i \psi(B)
\]
for any self-adjoint operators $A,B \in F(H)$
is a Jordan homomorphism of $F(H)$.
See, for example, the proof of \cite[Theorem 2]{MolStud1}.

Since $F(H)$ is a locally matrix ring, it follows from a celebrated
result
of Jacobson and Rickart \cite[Theorem 8]{JR} that $\tilde \psi$ can
be written as $\tilde \psi= \psi_1 +\psi_2$, where $\psi_1$ is a
homomorphism and $\psi_2$ is an antihomomorphism. Let $P$ be a rank-one
projection. Since $\tilde \psi(P)=\phi(P)$ is also rank-one, we obtain
that
one of the idempotents $\psi_1(P), \psi_2(P)$ is zero. Since $F(H)$ is a
simple ring, it is now easy to see that this implies that either
$\psi_1$ or $\psi_2$ is identically zero, that is, $\tilde \psi$ is
either a homomorphism or an antihomomorphism of $F(H)$. In what follows
we can assume without loss of generality that $\tilde \psi$ is a
homomorphism.

We show that $\tilde \psi$ preserves the rank. Let $A$ be a
rank-$n$ operator. Then there is a rank-$n$ projection $P$ such that
$PA=A$.
Thus, $\tilde \psi(A)=\tilde \psi(P)\tilde \psi(A)=\phi(P)\tilde
\psi(A)$
which proves that $\tilde \psi(A)$ is of rank at most $n$. If $Q$ is
any rank-$n$ projection, then there are finite rank operators $U,V$
such that $Q=UAV$. Since $\phi(Q)=\tilde \psi(Q)=\tilde \psi(U) \tilde
\psi (A)\tilde
\psi(V)$ and the rank of $\phi(Q)$ is $n$, it follows that the rank of
$\tilde \psi(A)$ is at least $n$. Therefore, $\tilde \psi$ is rank
preserving. We now refer to Hou's work \cite{Hou}. It follows from the
argument leading to \cite[Theorem 1.2]{Hou} (which is in fact a standard
'preserver-argument' already) that there are linear
operators $T,S$ on $H$ such that $\tilde \psi$ is of the form
\begin{equation}\label{E:gleasegy}
\tilde \psi(x\ot y)=(Tx)\ot (Sy) \qquad (x,y \in H)
\end{equation}
(recall that we have assumed that $\tilde \psi$ is a homomorphism).
We claim that $T,S$ are bounded. This will follow from the following
lemma.

\begin{lemma}\label{L:cont}
Let $T,S$ be linear operators on $H$ with the property that the map $x
\mapsto (Tx) \ot (Sx)$ is continuous on the unit ball of $H$.
Then $T, S$ are bounded.
\end{lemma}

\begin{proof}
If $x_n \to x$ and $x\neq 0$, then we have
\[
\frac{x_n}{\| x_n\|} \to \frac{x}{\| x\|}.
\]
This implies that
\[
\frac{(Tx_n) \ot (Sx_n)}{\| x_n\|^2} \to \frac{(Tx)\ot (Sx)}{\| x\|^2}
\]
which yields
\[
(Tx_n) \ot (Sx_n) \to (Tx)\ot (Sx).
\]
Consequently, the map $x \mapsto (Tx)\ot (Sx)$ is continuous at
any point different from 0.
Now, let $x_n \to 0$ and pick a nonzero vector $y\in H$ for which $Sy
\neq
0$ (observe that if $S=0$, then there is nothing to prove). Using the
polarization identity
\[
\begin{split}
(Tx_n)\ot (Sy)= &
{\textstyle \frac{1}{4}} \{ T(x_n+y)\ot S(x_n+y)-T(x_n-y)\ot
S(x_n-y)+\\ & iT(x_n+iy)\ot S(x_n+iy)-iT(x_n-iy)\ot S(x_n-iy)\} ,
\end{split}
\]
we see that $(Tx_n) \ot (Sy) \to 0$ which gives us that $T$ is
continuous
at 0, that is, $T$ is bounded. The boundedness of $S$ is now obvious.
\end{proof}

To continue the proof of Theorem~\ref{T:multrankPH}, we infer from
\eqref{E:gleasegy} that $\la Tx, Sx\ra=\la x,x\ra$ for every unit vector
$x\in H$ ($\phi$ sends rank-one projections to idempotents). Clearly,
this implies that $\la Tx, Sy\ra=\la x,y\ra$ $(x,y\in H)$.
Consequently $S^*T=I$. We have $\phi(P)=TPS^*$ for every
rank-one
projection $P\in P(H)$. By the additivity property of $\phi$ appearing
in \eqref{E:hat}, it follows that $\phi(P)=TPS^*$
holds true for every finite-rank projection $P$ as well.

Denote $Q=TS^*$. Clearly, $Q^2=TS^*TS^*=TIS^*=Q$. Let $P$ be an
arbitrary projection. Choose a monotone increasing sequence $(P_n)$ of
finite-rank projections which weakly
converges to $I$. We compute
\begin{equation}\label{E:gleasket}
\phi(P)Q= \phi(P) \wlim_n TP_nS^*= \wlim_n \phi(P)TP_nS^*=
\end{equation}
\[
\wlim_n \phi(P)\phi(P_n)= \wlim_n \tilde \psi(PP_n)= \wlim_n TPP_nS^*=
TPS^*
\]
and
\begin{equation}\label{E:gleasketa}
Q \phi(P)= (\wlim_n TP_nS^*)\phi(P)= \wlim_n (TP_nS^*\phi(P))=
\end{equation}
\[
\wlim_n \phi(P_n)\phi(P)= \wlim_n \tilde \psi(P_nP)= \wlim_n TP_nPS^*=
TPS^*.
\]
So, $Q$ is an idempotent commuting with the range of $\phi$. Therefore,
$\phi$ can be written as
\[
\phi(P)=\phi(P)Q + \phi(P)(I-Q)
\]
where the maps $\phi_1: P\mapsto \phi(P)Q$ and $\phi_2: P\mapsto
\phi(P)(I-Q)$ are multiplicative. We see that $\phi_1(P)=TPS^*$
$(P\in P(H))$ and thus $\phi_2$ vanishes on the set of all finite-rank
projections. This completes the proof of the theorem.
\end{proof}

\begin{proof}[Proof of Theorem~\ref{T:multrank}]
If we consider $\phi$ only on $P(H)$, then by Theorem~\ref{T:multrankPH}
it follows that
\begin{equation}
\phi(P)=TPS
\end{equation}
for every finite rank projection $P$,
where $T,S$ are either both bounded linear operators or both bounded
conjugate-linear operators with $ST=I$. In what follows
we can suppose without loss of generality that $T,S$ are linear.

Let $A\in B(H)$ be a rank-one operator. Then there is another rank-one
operator
$R$ such that $A=ARA$. Since $\phi$ preserves the rank, it follows from
the equality
\[
\phi(\l A)=\phi(A) \phi(\l R) \phi(A)
\]
that $\phi(\l A)=f_A(\l) \phi(A)$ with some scalar function $f_A$.
If $B$ is a rank-one operator with $BA\neq 0$ and $f_B$ is
the scalar function corresponding to $B$, then we have
\[
f_A(\l)\phi(B)\phi(A)=
\phi(B) \phi(\l A)=
\phi(\l B)\phi(A)=f_B(\l) \phi(B)\phi(A)
\]
which implies that $f_A(\l)=f_B(\l)$ $(\l \in \C)$.
If $C$ is a rank-one operator and $CA=0$,
then we can choose a rank-one operator $B$ such that $CB\neq 0$ and
$BA\neq 0$. This gives us that $f_C=f_B=f_A$. Therefore, the scalar
function $f_A$ does not depend on the rank-one operator $A$. In what
follows this function will be denoted by $f$.
It follows from the equality
\[
f(\l \mu)\phi(A)=
\phi(\l \mu A)=f(\l)\phi(\mu A)=f(\l) f(\mu) \phi(A)
\]
that $f$ is a continuous multiplicative function. We show that it is
additive as well.
Let $x, y\in H$ be orthogonal unit vectors. Since $\phi$ is additive
on the set of finite rank projections, we compute
\[
\phi((\l x + \mu y)\ot y)=
\phi( x\ot x +y\ot y)\phi((\l x + \mu y)\ot y)=
\]
\[
(\phi( x\ot x )+\phi(y\ot y))\phi((\l x + \mu y)\ot y)=
\]
\[
\phi( x\ot x )\phi((\l x + \mu y)\ot y)+
\phi(y\ot y)\phi((\l x + \mu y)\ot y)=
\]
\[
\phi(\l x \ot y)+\phi(\mu y\ot y)=
f(\l)\phi( x \ot y)+f(\mu)\phi(y\ot y)
\]
Multiplying by $\phi(x \ot (x+y))$ from the left we can compute
\[
f(\l+\mu) \phi(x\ot y)=
\phi((\l+\mu)(x\ot y))=
\]
\[
\phi((x\ot (x+y))(\l x + \mu y)\ot y)=
\phi(x\ot (x+y))\phi((\l x + \mu y)\ot y)=
\]
\[
\phi(x\ot (x+y))
(f(\l)\phi( x \ot y)+f(\mu)\phi(y\ot y))=
\]
\[
f(\l)\phi( (x\ot (x+y))(x \ot y))+f(\mu)\phi((x\ot (x+y))(y\ot y))=
\]
\[
f(\l)\phi( x \ot y)+f(\mu)\phi(x\ot y).
\]
It follows that $f(\l+\mu)=f(\l)+f(\mu)$, that is, $f$ is additive.
Therefore, $f$ is a continuous ring endomorphism of $\C$ with $f(1)=1$.
It is well-known that this implies that $f$ is either the identity or
the conjugation on $\C$. We show that in our case $f$ is the identity.
Suppose on the contrary that $f(\l )=\bar \l$ $(\l \in \C)$.
Let $x,y$ be non-orthogonal unit vectors. Since $ST=I$, we compute
\[
\la x,y \ra \phi(x\ot y)=
\phi(\la y,x \ra x\ot y)=
\phi(x \ot x \cdot y\ot y)=
\]
\[
\phi(x\ot x)\phi(y \ot y)=Tx\ot S^*x \cdot Ty \ot S^*y=
\]
\[
\la Ty,S^*x\ra Tx \ot S^*y=
\la y,x\ra Tx \ot S^*y.
\]
So, we have
\begin{equation}\label{E:tiz}
\phi(x\ot y)= \frac{\la y,x\ra}{\la x,y \ra} Tx \ot S^*y.
\end{equation}
Now, let $x,y,u,v$ be unit vectors for which $\la x,y \ra, \la x,v \ra,
\la u,v\ra \neq 0$. We then have
\[
\phi(x \ot y )\phi(u\ot v)=
\frac{\la y,x\ra}{\la x,y \ra} Tx \ot S^*y \cdot
\frac{\la v,u\ra}{\la u,v \ra} Tu \ot S^*v=
\]
\[
\frac{\la y,x\ra}{\la x,y \ra}\frac{\la v,u\ra}{\la u,v \ra}
\la u,y\ra  Tx \ot S^*v.
\]
On the other hand,
\[
\phi(x \ot y \cdot u\ot v)=\phi(\la u,y\ra x\ot v)=
\la y,u\ra \phi(x\ot v)=
\la y,u\ra \frac{\la v,x\ra}{\la x,v \ra} Tx \ot S^*v.
\]
Comparing these two equalities we arrive at
\[
\la y,x\ra  \la u,y\ra \la v,u\ra \la x,v \ra=
\la x,y \ra \la y,u\ra \la u,v\ra \la v,x\ra.
\]
Since this equality obviously does not hold true for every possible
choice of $x,y,u,v\in H$, we obtain that $f$ is really the identity.

Now, the same argument that has led to \eqref{E:tiz}
shows that
\[
\phi(x\ot y)= Tx \ot S^*y
\]
if $x,y\in H$ are non-orthogonal unit vectors. If $x,y$ are orthogonal,
then
choosing a unit vector $z\in H$ such that $\la x,z\ra\neq 0$ and $\la
z,y\ra \neq 0$ we have
\[
\phi(x\ot y)=\phi(x\ot z)\phi(z\ot y)=
Tx \ot S^*z \cdot Tz \ot S^*y=
\]
\[
\la Tz,S^*z \ra Tx \ot S^*y =
\la z,z \ra Tx \ot S^*y =
Tx \ot S^*y.
\]
Since $f$ is the identity, we thus obtain
$\phi(A)=TAS$ for every rank-one operator $A$.
If $A\in F(H)$ and $P$ is a finite-rank projection such that $A=PA$
and $P_1, \ldots, P_n$ are pairwise orthogonal rank-one projections
such that $P=P_1 +\cdots +P_n$, then it follows that
\[
\phi(A)=\phi(PA)=\phi(P)\phi(A)=
\]
\[
\sum_i \phi(P_i) \phi(A)=
\sum_i \phi(P_i A)=
\sum_i T(P_iA)S= TAS.
\]
Similarly as in \eqref{E:gleasket}, \eqref{E:gleasketa} in the proof of
Theorem~\ref{T:multrankPH} we see that
the operator $Q=TS$ is an idempotent commuting with the range of $\phi$
and $\phi(A)Q=TAS$ $(A\in B(H))$. Therefore, $\phi$ can be written as
\[
\phi(A)=\phi(A)Q + \phi(A)(I-Q)
\]
where the maps $\phi_1: A\mapsto \phi(A)Q$ and
$\phi_2: A\mapsto \phi(A)(I-Q)$ are multiplicative and $\phi_2$ vanishes
on the set of all finite rank operators.
We claim that $\phi_2$ is identically 0.
Indeed, if $\phi_2$ is not zero, then $\phi_2(I)\neq 0$. If $P$ is a
projection of infinite rank, then due to the fact that in that case
there is a coisometry $U$ such that $UPU^*=I$, it follows that
$\phi_2(P)\neq 0$. Choosing an uncountable set of infinite rank
projections in $B(H)$ with the property that the product of any two of
them has finite rank (see the first part of the proof of \cite[Theorem
1]{MolStud1}) and taking the values of those
projections under $\phi_2$, we would obtain uncountably many pairwise
orthogonal nonzero idempotents in $B(H)$ which contradicts the
separability of $H$. This shows that $\phi_2=0$.
So, $\phi(A)=\phi(A)Q=TAS$ for every $A\in B(H)$.
This completes the proof.
\end{proof}

\begin{proof}[Proof of Theorem~\ref{T:multcorank}]
We prove that $\phi$ preserves the rank of projections. First suppose
that $\phi(P)=0$ for every finite rank projection $P\in B(H)$. Since
$\phi(I)=I$, just as in the proof of Theorem~\ref{T:multrank} we see
that $\phi(P)\neq 0$ for every infinite rank projection $P$ and then we
arrive at a contradiction in the same way as there. So, let $n$ be the
smallest positive integer with the property that $\phi(P)\neq 0$
whenever $P\in P(H)$ is of rank $n$ (observe that by the
multiplicativity of $\phi$, we have
$\rank \phi(Q)= \rank \phi(Q')$ if $\rank Q=\rank Q'$).
We claim that the rank of $\phi(P)$ is 1 for every such $P$.
Indeed, let
$Q$ be a rank-one projection and $P$ be a rank-$n$ projection such that
$(I-Q)P=P(I-Q)$ is of rank $n-1$. Then
$\phi(I-Q)$ and $\phi(P)$ are orthogonal and we have
$\phi(I-Q) +\phi(P)\leq I$. Since the corank of $\phi(I-Q)$ is 1, this
gives us that the rank of $\phi(P)$ is
1. We show that $\rank P=1$. Suppose on the contrary that
$\rank P=n>1$. Let $P\leq R$ be a projection of rank $n+1$.
Similarly as just
before, we can verify that the rank of $\phi(R)$ is at most 2. On the
other hand there are rank-$n$ projections $P_1, \ldots, P_{n+1}\leq R$
such that the product of any two of them is a rank-$(n-1)$ projection.
Consequently, $\phi(P_1), \ldots,
\phi(P_{n+1})$ are orthogonal and $\phi(P_1)+\cdots +\phi(P_{n+1})\leq
\phi(R)$. Therefore, we have $n+1 \leq 2$. This gives us that
$n=1$ and hence $\phi$ sends rank-one projections to rank-one
idempotents.

Let now $P$ be a rank-$n$
projection. Since $\phi(I-P)$, $\phi(P)$ are orthogonal
idempotents and $\phi(I-P)$ has corank $n$, we obtain that $\phi(P)$ has
rank at most $n$. Now, if $P_1, \ldots, P_n$ are pairwise orthogonal
rank-1 projections, then like in the proof of Theorem~\ref{T:multrankPH}
we see that
\[
\phi(P_1)+\cdots +\phi(P_n)\leq
\phi(P_1+\cdots +P_n).
\]
Since the idempotent appearing on the left-hand side of this inequality
has
rank $n$ and the one on the right-hand side has rank at most $n$, we
infer that
\[
\phi(P_1)+\cdots +\phi(P_n)=
\phi(P_1+\cdots +P_n).
\]
Therefore, $\phi$ preserves the rank of projections. Similarly to the
argument in the proof of Theorem~\ref{T:multrankPH} before
Lemma~\ref{L:cont} we get that $\phi$ is rank-preserving. By
Theorem~\ref{T:multrank} we have the form
\eqref{E:formhar} of $\phi$. Since $\phi(I)=I$, we also obtain $TS=I$.
\end{proof}

\begin{proof}[Proof of Theorem~\ref{T:multspect}]
Let $\l \in \C$. By the properties of $\phi$, $\phi(\l I)$ is a normal
operator whose spectrum does not contain any scalar different from $\l$.
Therefore, $\phi(\l I)=\l I$. This gives us that $\phi$ is homogenous.

We prove that for any orthogonal projections $P, Q$ we have
$\phi(P+Q)=\phi(P)+\phi(Q)$. Let $P, Q$ be of infinite rank such that
$P+Q=I$. Pick a scalar $0< \mu <1$. We have
$\sigma (\phi(P+\mu Q)) \subset \{ 1, \mu\}$.
We distinguish three cases.
First suppose that
$\sigma (\phi(P+\mu Q)) = \{ 1\}$. Since $\phi$ preserves normality,
this yields $\phi(P+\mu Q)=I$. Taking powers, we obtain
\[
\phi( P+\mu ^n Q)=(\phi(P+\mu Q))^n =I \qquad (n\in \N).
\]
Using the continuity of $\phi$ we have $\phi(P)=I$.
Since $\phi(P)+\phi(Q) \leq I$, we get $\phi(Q)=0$. On the other
hand,
$P, Q$ are equivalent projections and it follows that $\phi(P)=0$
which is a contradiction.
Next suppose that $\sigma (\phi(P+\mu Q)) = \{ \mu \}$. Then we have
$\phi(P+\mu Q)=\mu I$. Taking powers again, we have
\[
\phi( P+\mu ^n Q)=(\phi(P+\mu Q))^n =\mu^n I \longrightarrow 0 \qquad
(n\in \N).
\]
Thus, we infer $\phi(P)=0$ which gives us that $\phi(I)=0$, a
contradiction.
Consequently, we have $\sigma (\phi(P+\mu Q)) = \{ 1, \mu \}$.
This implies that $\phi (P+\mu Q)= P'+\mu Q'$ where $P',Q'$ are nonzero
projections such that $P'+Q'=I$. We show that in this case
$\phi (P+\epsilon Q)= P'+\epsilon Q'$ holds for every $0< \epsilon <1$.
By what we have just proved, for an arbitrary $0<\epsilon<1$ we can
write $\phi (P+\epsilon Q)$ as
\[
\phi (P+\epsilon Q)= P''+\epsilon Q''
\]
where $P'', Q''$ are orthogonal nonzero projections with $P''+Q''=I$.
Since $\phi$ clearly preserves the commutativity, we get that
$P'+\mu Q'$ and $P''+\epsilon Q''$ commute. Referring to the spectral
theorem we obtain that $P', Q', P'', Q''$ are pairwise commuting.
Furthermore, since
\[
\phi (P+\mu Q)\phi (P+\epsilon Q)=\phi (P+\mu \epsilon Q),
\]
it follows that
$(P'+\mu Q')(P''+\epsilon Q'')$ is of the form
$P'''+\mu \epsilon Q'''$. Because of the equality
\[
(P'+\mu Q')(P''+\epsilon Q'')=P'P''+\epsilon P'Q'' +\mu Q'P'' +\epsilon
\mu Q'Q''
\]
and the fact that the spectrum of $\phi (P+\mu \epsilon Q)$ is $\{ 1,
\mu\epsilon\}$ we obtain that $P'Q''=Q'P''=0$. Therefore, $P'\leq P''$
and $Q'\leq Q''$. Since $P'+Q'=P''+Q''=I$, it follows that $P'=P''$ and
$Q'=Q''$. So, we have
\[
\phi (P+\mu Q)= P'+\mu Q'
\]
for every $0< \mu <1$. Sending $\mu$ to 0, we get
$\phi(P)=P'$. Since $\phi(I)=I$, $\phi$ preserves the inverse operation.
This yields that
\[
\phi (P+(1/\mu) Q)= P'+(1/\mu) Q'.
\]
By the homogenity of $\phi$ we infer that
\[
\phi (\mu P+ Q)= \mu P'+Q'.
\]
If $\mu \to 0$, we arrive at $\phi(Q)=Q'$. Consequently, we have
\begin{equation}\label{E:nyolc}
\phi(P)+\phi(Q)=I.
\end{equation}
If $R,R'$ are projection such that $R\leq P$
and $R'\leq Q$, then multiplying \eqref{E:nyolc} by $\phi(R+R')$ we
arrive at
\begin{equation}\label{E:kilenc}
\phi(R)+\phi(R')=\phi(R+R').
\end{equation}
Therefore, we have \eqref{E:kilenc} whenever
$R,R'$ are orthogonal, either both infinite or both finite rank
projections. If $R$ is of finite rank and
$R'$ is of infinite rank, then we can write $R'=P+Q$ where $P,Q$ are
orthogonal and they are of infinite rank. The argument leading
to \eqref{E:kilenc} gives us that $\phi(R')=\phi(P)+\phi(Q)$ and
$\phi(R+P)=\phi(R)+\phi(P)$. We then have
\[
\phi(R+R')= \phi(R+P+Q)=\phi(R+P)+\phi(Q)=
\]
\[
\phi(R)+\phi(P)+\phi(Q)= \phi(R)+\phi(R').
\]
Hence, $\phi$ is additive on the set of all projections. Since
$\phi$ sends projections to projections, $\phi$ is bounded on $P(H)$. By
a deep result due to Bunce and Wright \cite[Theorem A]{BW} it follows
that $\phi_{| P(H)}$ can be extended to a bounded linear transformation
$\psi$ on $B(H)$. Since $\psi$ sends projections to projections and
$\psi$ is
continuous, it is a standard argument to verify that $\psi$ is a Jordan
*-homomorphism (once again, see the proof of \cite[Theorem
2]{MolStud1}).

We next refer to the proof of \cite[Theorem 3]{BaMo}.
Similarly to the argument followed there,
we obtain that there is a central projection $Q$ in the
$C^*$-algebra generated by the range
of $\psi$ such that $\psi_1(.)=\psi(.)Q$ is a *-homomorphism and
$\psi_2(.)=\psi(.)(I-Q)$ is a *-antihomomorphism.
This gives us that $\psi$ is the direct sum of the maps
\begin{equation}\label{E:tizen}
\psi_1: A \longmapsto \sum_n U_n A U_n^*
\end{equation}
and
\[
\psi_2: A \longmapsto \sum_n V_n A^{tr} V_n^*,
\]
where $U_n,V_n:H\to H$ are isometries with pairwise orthogonal ranges
and $^{tr}$ denotes the transpose with respect to a fixed orthonormal
basis in $H$. Consequently,
$\psi$ can be represented as \begin{equation}\label{E:tizenegy}
\psi(A)=
\left[
\begin{matrix}
A       &   0      & \dots     & 0        & 0       &\dots \\
0       &   A      & \dots     & 0        & 0       &\dots \\
\vdots  & \vdots   & \ddots    & \vdots   & \vdots  &\vdots\\
0       &   0      & \dots     & A^{tr}   & 0       &\dots \\
0       &   0      & \dots     & 0        & A^{tr}  &\dots \\
\vdots  & \vdots   & \vdots    & \vdots   & \vdots  &\ddots \\
\end{matrix}
\right].
\end{equation}
We show that the *-antihomomorphic part of $\psi$ is in fact missing,
and hence $\psi$ is a *-homomorphism. Let $P_1, \ldots, P_n$ be
pairwise orthogonal projections and let $\l_1, \ldots, \l_n$ be
scalars. We compute
\[
\phi(\sum_i \l_i P_i)=
\phi((\sum_i \l_i P_i)(\sum_i P_i))=
\phi(\sum_i \l_i P_i)\phi(\sum_i P_i)=
\]
\[
\phi(\sum_i \l_i P_i)\sum_i \phi(P_i)=
\sum_k \phi(\sum_i \l_i P_i) \phi(P_k)=
\sum_k \phi((\sum_i \l_i P_i) P_k)=
\]
\[
\sum_k \phi(\l_k P_k)=
\sum_k \l_k \phi(P_k)=
\psi(\sum_k \l_k P_k).
\]
By the continuity of $\phi, \psi$ and the spectral theorem we get that
$\phi(N)=\psi(N)$ holds for every normal operator $N\in B(H)$.
Suppose that $^{tr}$ do appear in \eqref{E:tizenegy}.
If $S_1, \ldots, S_n$ are self-adjoint operators such that $N=S_1 \cdot
\ldots \cdot S_n$ is normal, then by the multiplicativity of $\phi$ we
have \[
N^{tr}=S_1^{tr} \cdot \ldots \cdot S_n^{tr}=(S_n \cdot \ldots \cdot
S_1)^{tr}={N^*}^{tr}
\]
which yields that $N=N^*$. Let $x, y$ be orthogonal unit vectors and
$S_1=x\ot y+y \ot x$, $S_2=x\ot x -y\ot y$. It is trivial to check that
$N=S_1S_2$ is normal but not self-adjoint.
Therefore, we obtain that $\psi_2=0$, that is, $\psi$ is a
*-homomorphism.
It is easy to see that every rank-one operator is the scalar multiple of
the product of at most three rank-one projections. This gives us that
$\psi$ and $\phi$ coincide on the rank-one operators. To complete the
proof, let $A\in B(H)$ be arbitrary. Choose a maximal set $(P_n)$ of
pairwise orthogonal rank-one projections in $B(H)$. We compute
\[
\phi(A)=\phi(A)\phi(I)=\phi(A)\psi(I)=\phi(A)\sum_n \psi(P_n)=
\phi(A)\sum_n \phi(P_n)=
\]
\[
\sum_n \phi(A)\phi(P_n)=
\sum_n \phi(AP_n)=
\sum_n \psi(AP_n)=
\psi(A),
\]
where we have used the weak continuity of $\psi$ which clearly holds by
\eqref{E:tizen}. Finally, since $\phi(I)=I$, we have $\sum_n
U_nU_n^*=I$. This completes the proof.
\end{proof}

\bibliographystyle{amsplain}

\end{document}